\documentclass[12pt]{amsart}

\usepackage{amsmath,amsthm,amssymb}
\font\teneufm=eufm10
\font\seveneufm=eufm7
\font\fiveeufm=eufm5
\newfam\frakturfam

\textfont\frakturfam=\teneufm
\scriptfont\frakturfam=\seveneufm
\scriptscriptfont\frakturfam=\fiveeufm

\newtheorem{pr}{Proposition}

\newtheorem{lm}{Lemma}
\newtheorem{theor}{Theorem}
\newtheorem{co}{Corollary}

\def\bee{\begin{eqnarray}}
\def\bes{\begin{eqnarray*}}
\def\eee{\end{eqnarray}}
\def\ees{\end{eqnarray*}}

\def\a{\alpha}
\def\b{\beta}

\def\Proof{{\sl Proof.}\ }

\title{Representations via differential algebras and equationally Noetherian algebras}

\begin{document}
\date{}
\maketitle

\begin{center}
{\bf Alexander A.~Mikhalev}\footnote{Faculty of Mechanics
and Mathematics, 
M.V.Lomonosov 
Moscow State University, Leninskie Gory~1,
119991, Moscow, Russian Federation 
e-mail: {\em aamikhalev@mail.ru}.
},
{\bf Manat Mustafa}\footnote{Department of Mathematics, SSH., Nazarbayev University, 
53 Kabanbay Batyr Avenue, Astana, 010000, Kazakhstan,
e-mail: {\em manat.mustafa@nu.edu.kz}},
and
{\bf Ualbai Umirbaev}\footnote{Department of Mathematics,
 Wayne State University,
Detroit, MI 48202, USA; 
and Institute of Mathematics and Mathematical Modeling, Almaty, 050010, Kazakhstan,
e-mail: {\em umirbaev@wayne.edu}}

\bigskip

\dedicatory{To the memory of A.~V.~Mikhalev (1940--2022)}

\end{center}

\begin{abstract} We show that free algebras of the variety of algebras generated by the Witt algebra $W_n$, the left-symmetric Witt algebra $L_n$, and the symplectic Poisson algebra $P_n$  can be described as subalgebras of differential polynomial algebras with respect to  appropriately defined products. Using these representations, we prove that $W_n$, $L_n$, $P_n$, and  the free algebras of the varieties of algebras generated by these algebras are equationally Noetherian.
\end{abstract}

\noindent


\section{Introduction}

\hspace*{\parindent}

The well-known Hilbert Basis Theorem (see, for example \cite{Eisenbud}) implies that the affine space $\mathbb{A}^n(K)$ over an arbitrary field $K$ is Noetherian. In other words, any field $K$ is {\em equationally Noetherian}.  The basics of algebraic
geometry over groups and algebraic systems can be found in \cite{BMR,DMR12,Plotkin}. Many standard notions of algebraic geometry such as algebraic
sets, the Zariski topology, irreducible
varieties, and coordinate groups were freely used to
solve the famous Tarski's problems independently by O. Kharlampovich and A. Myasnikov  \cite{KM05,KM06} and by Z. Sela \cite{Sela06}. In 1986 V. Guba \cite{Guba} proved that free groups are equationally Noetherian using the well known representation of free groups via matrices of order two. It is also known that  linear groups over a Noetherian domain \cite{BMR} (in particular, polycyclic groups and finitely generated metabelian groups), finitely generated free solvable groups \cite{GR}, and torsion-free hyperbolic groups \cite{Sela09} are equationally Noetherian.

There are only few examples of equationally Noetherian rings and algebras. Obviously, Hilbert's basis theorem implies that every associative and commutative Noetherian ring or algebra is equationally Noetherian.
It is also known
that any finitely generated metabelian or nilpotent Lie algebra is 
equationally Noetherian \cite{DMR12}. Finitely generated free metabelian Lie algebras of characteristic zero  are Noetherian with respect to equations in the universal enveloping algebras \cite{RS}.

In this paper we give more examples of equationally Noetherian algebras. The source of these examples is differential algebras.  The Ritt-Raudenbush Basis Theorem \cite{Kaplansky, Kolchin, Ritt}, which is an analogue of the Hilbert Basis Theorem, says that every finitely generated differential algebra of characteristic zero is radically Noetherian, i.e., satisfies the ascending chain condition for radically closed differential ideals. Consequently, every  finitely generated differential domain over a field of characteristic zero is equationally Noetherian. 

We show that free algebras of the variety of algebras generated by the Witt algebra $W_n$, the left-symmetric Witt algebra $L_n$, and the symplectic Poisson algebra $P_n$  over a field of characteristic zero can be described as subalgebras of differential polynomial algebras with respect to  appropriately defined products.

Using these representations, we prove that the algebras $W_n$, $L_n$, $P_n$, and the free algebras of the varieties of algebras generated by these algebras are equationally Noetherian.

Recall that the variety of algebras generated by the left symmetric Witt algebra $L_1$ is the variety of all Novikov algebras \cite{MLU11N}. Consequently, every free Novikov algebra over a filed of characteristic zero is equationally Noetherian.  The solvability, nilpotency, and identities of Novikov algebras have been intensively studied in recent years  \cite{ShZh,TUZ21,UZh21,ZhU}. In particular, every Novikov algebra over a field of characteristic zero satisfying a nontrivial polynomial identity is associative right nilpotent \cite{DIU}.

The paper is organized as follows. In Section 2 we give an exact definition of equationally Noetherian algebras, formulate some basic results from \cite{DMR12}, and give some examples of equationally Noetherian free algebras. In Section 3 we recall some well known results on differential algebras from \cite{Kaplansky, Kolchin, Ritt} and give some examples of equationally Noetherian differential algebras.  
Section 4 is devoted to a representation of free algebras of the variety of algebras generated by the left-symmetric Witt algebra $L_n$ and the Witt algebra $W_n$. Here we also prove that  $W_n$, $L_n$, and the free algebras of the varieties of algebras generated by these algebras are equationally Noetherian. The same results for the symplectic  Poisson algebra $P_n$ are formulated in Section 5.

\section{On equationally Noetherian algebras}

\hspace*{\parindent}

 Let $\mathfrak{M}$ be an arbitrary variety of (universal) algebras  and let 
$A$ be an arbitrary  algebra  of $\mathfrak{M}$. Let $\mathfrak{M}\langle x_1,x_2,\ldots,x_m\rangle$ be the free algebra of the variety $\mathfrak{M}$ freely generated by $x_1,x_2,\ldots,x_m$.  Consider the coproduct (or free product)
\bes
B=A\sqcup_{\mathfrak{M}}\mathfrak{M}\langle x_1,x_2,\ldots,x_m\rangle
\ees
 of $A$ and $\mathfrak{M}\langle x_1,x_2,\ldots,x_m\rangle$ in the variety of algebras $\mathfrak{M}$. Recall that (see, for example \cite{Mac}) there exist homomorphisms $\a_1 : A\to B$ and $\a_2 : \mathfrak{M}\langle x_1,x_2,\ldots,x_m\rangle\to B$ satisfying the following universal property: for any algebra $C\in \mathfrak{M}$ and for any two homomorphisms $\phi_1 : A\to C$ and $\phi_2 : \mathfrak{M}\to C$ there exists a unique homomorphism $\phi : B\to C$ such that $\phi \circ \a_1=\phi_1$ and 
$\phi \circ \a_2=\phi_2$. Notice that the homomorphisms $\a_1$ and $\a_2$ are not always injections. In this paper we consider only the varieties of algebras for which the homomorphisms $\a_1$ and $\a_2$ are always injections. Then we can identify $A$ 
and $\mathfrak{M}\langle x_1,x_2,\ldots,x_m\rangle$ with there images in $B$. As a consequence of this, for any $m$-tuple $(a_1,\ldots,a_m)\in A^m$ there exists the unique homomorphism $\phi : B\to A$ such that $\phi_A=\mbox{id}$ and $\phi(x_i)=a_i$ for all $i$. For any $f=f(x_1,\ldots,x_m)\in B$ we set $f(a_1,\ldots,a_m)=\phi(f)$. 
This allows to consider any element $f=f(x_1,\ldots,x_m)$ of $B$ as a function 
\bes
f : A^m\to A. 
\ees

 For any subset $S$ of $B$ we define the set of all {\em zeroes} 
\bes Z(S)=\{(a_1,\ldots,a_m) | f(a_1,\ldots,a_m)=0, f\in
S\} 
\ees of $S$ in $A^m$. The sets $Z(S)$ allow to define the Zariski
topology on $A^m$. The algebra $A$ is called {\em equationally Noetherian} if $A^m$ is Noetherian with respect to Zariski
topology for all natural $m$.

Notice that this definition
does not depend on the choice of $\mathfrak{M}$ such that $A\in\mathfrak{M}$. In many cases it is more efficient to assume that $\mathfrak{M}$ is the variety of algebras $\mathbf{Var}(A)$ generated by $A$. Sometimes it is better to choose $\mathfrak{M}$ with some special properties, for example, the varieties with an explicit description of the coproduct of algebras. It might happen that $A\sqcup B=0$ for nonzero algebras $A$ and $B$ in some varieties of algebras with unity in the signature.

\begin{pr}\label{p1} \cite{DMR12} Let $A$ be an equationally Noetherian algebra. Then the following algebras (of the same signature) are equationally Noetherian:

(1) Every subalgebra of $A$;

(2) Every direct power of $A$;

(3) Every algebra of the quasi variety of algebras $\mathbf{Qvar}(A)$ generated by $A$;

(4) Every free algebra of the variety of algebras $\mathbf{Var}(A)$ generated by $A$.
\end{pr}

Let us revise some examples of equationally Noetherian associative and commutative algebras:

(i) Every field is equationally Noetherian. This is a corollary of Hilbert's basis theorem and is one of the basics of classical algebraic geometry.

(ii) Every domain is equationally Noetherian by Proposition \ref{p1}(1) since it is embeddable into a field.

(iii) Every Noetherian ring with identity is also equationally Noetherian by Hilbert's basis theorem.

(iv) Let $R$ be an arbitrary equationally Noetherian associative and commutative algebra with identity containing the field of rationals $\mathbb{Q}=\mathbb{Q}\cdot 1\subseteq R$. Then every polynomial algebra over $R$ is equationally Noetherian. 
This follows from Proposition \ref{p1}(4) since $R$ generates the variety of all associative and commutative algebras over $R$. 

The next theorem gives more examples of associative and nonassociative equationally Noetherian algebras.

 \begin{theor}\label{t1} Let $R$ be an arbitrary associative and commutative equationally Noetherian ring with identity. Let $A$ be an arbitrary $R$-algebra that is a finitely generated free module over $R$. Then $A$ is an equationally Noetherian $R$-algebra.
\end{theor}
\Proof For simplicity of notations, we prove this theorem only for rings $R$ containing the field of rational numbers $\mathbb{Q}$. 

Let $e_1,\ldots,e_m$ be a basis of a finite dimensional algebra $A$ over $R$. Then every element $a\in A$ can be identified with $(r_1,\ldots,r_m)\in R^m$ such that
\bes
a=r_1e_1+\ldots+r_me_m, 
\ees
i.e., 
\bes
A\equiv R^m
\ees
and 
\bee\label{f1}
A^n\equiv R^{mn}.
\eee

Let $R_0=R[x_i^j : 1\leq i\leq n, 1\leq j\leq m]$ be a polynomial algebra over $R$. Consider the algebra $A_0=R_0\otimes_{R}A$. Since $A_0$ is a tensor product of free modules, $R_0$ and $A$ are embeddable into $A_0$. 
Let $\mathfrak{M}=\mathrm{Var}(A)$ be the variety of algebras over $R$ generated by $A$. Since $\mathbb{Q}\subseteq R$, it follows that $A_0\in \mathfrak{M}$. 

Set 
\bes
X_i=x_i^1\otimes e_1+x_i^2\otimes e_2+\ldots +x_i^m\otimes e_m\in A_0, \ \ \ 1\leq i\leq n.
\ees
Let $R\langle X_1,\ldots,X_n\rangle$ be the $R$-subalgebra of $A_0$ generated by 
$X_1,\ldots, X_n$. It is well known that $R\langle X_1,\ldots,X_n\rangle$ is the free algebra of the variety of algebras $\mathfrak{M}=\mathrm{Var}(A)$ freely generated by $X_1,\ldots, X_n$. 

Let $C$ be the $R$-subalgebra of $A_0$ generated by $A$ and 
$R\langle X_1,\ldots,X_n\rangle$. Denote by $D=A\sqcup R\langle X_1,\ldots,X_n\rangle$ the coproduct of $A$ and $R\langle X_1,\ldots,X_n\rangle$ in $\mathfrak{M}$. Consider the homomorphism 
\bes
\omega : D\longrightarrow C
\ees
of $R$-algebras that is identical on $A$ and $R\langle X_1,\ldots,X_n\rangle$.

For any homomorphism $\varphi : D\rightarrow A$ with $\varphi_{|A}=\mathrm{id}$ there exists the unique homomorphism $\varepsilon: C\rightarrow A$ such that $\varepsilon_{|A}=\mathrm{id}$ and $\varphi=\varepsilon\circ \omega$. In fact, if
$\varphi(X_i)=a_i\in A$ and 
\bes
a_i=\sum r_i^j\otimes e_j, 
\ees
then we can define a homomorphism $\delta : A_0\to A$ that it is identical on $A$ and $\delta(x_i^j)=r_i^j$ for all $i,j$. We set $\varepsilon=\delta_{|C} : C\to A$. Then 
$\varepsilon(X_i)=a_i$ and $\varepsilon$ satisfies all required conditions. The uniqueness of $\varepsilon$ is obvious.

This means that if $f=f(X_1,\ldots,X_n)\in \mathrm{Ker}(\omega)$ then $f(a_1,\ldots,a_n)=0$ for all $a_1,\ldots,a_n\in A$. Consequently, the kernel of $\omega$ is not essential in studying of equations over $A$. Therefore it is sufficient to study equations of the form $s=0$ over $A$ where $s\in C$.

Let $S$ be an arbitrary set of elements of $C$. Every element  $s\in S$ can be uniquely written in the form
\bes
s=s_1\otimes e_1+\ldots+s_m\otimes e_m.
\ees
Using the identification (\ref{f1}), the equation $s=0$ over $A$ is equivalent to the equations $s_1=\ldots=s_m=0$ over $R$. Consequently,
the set of equations $s=0, s\in S,$ over $A$ is equivalent to the set of equations $s_1=\ldots=s_m=0, s\in S$ over $R$. The latter is equivalent to a finite subsystem of equations since $R$ is equationally Noetherian. Therefore the set of equations $s=0, s\in S,$ over $A$ is also equivalent to a finite subsystem, i.e., $A$ is equationally Noetherian. $\Box$

We formulate some essential corollaries of Theorem \ref{t1}.

 \begin{co}\label{c1} Every algebra of generic matrices over an equationally Noetherian associative and commutative algebra is equationally Noetherian.
\end{co}
\Proof Let $R$ be an equationally Noetherian associative and commutative ring. Then the algebra of matrices $M_n(R)$ of order $n$ over $R$ is equationally Noetherian by Theorem \ref{t1}. 
Recall that any algebra $G$ of generic matrices of order $n$ is a free algebra of the variety
of algebras $\mathrm{Var}(M_n(R))$. Then $G$ is equationally Noetherian by Proposition \ref{p1}. $\Box$

 \begin{co}\label{c2} Every free algebra of the variety of associative algebras defined by the identity
 \bes
 [x_1,y_1]\ldots [x_n,y_n]=0
 \ees
 over an equationally Noetherian associative and commutative algebra is equationally Noetherian.
\end{co}
\Proof It is well known that this variety of associative algebras over $R$ is generated by the algebra of triangular matrices $\mathrm{Tr}_n(R)$ \cite{Maltsev} since $\mathbb{Q}\subseteq R$. $\Box$

 \begin{co}\label{c3} Every free $\mathfrak{N}_c\mathfrak{A}$-Lie algebra
 over an infinite field is equationally Noetherian.
\end{co}
\Proof It is well known that this variety of Lie algebras is also generated by the algebra of triangular matrices $\mathrm{Tr}_n(R)$ \cite{Bahturin}. $\Box$

\section{Differential polynomial algebras}

\hspace*{\parindent}

The basic concepts of differential algebras can be found in \cite{Kaplansky,Kolchin,Ritt}.

Let $\Delta = \{\delta_1,\ldots, \delta_m\}$ be a basic set of derivation operators.

A commutative ring $R$ with identity is said to be a {\em differential} ring or $\Delta$-ring if all elements of $\Delta$ act on $R$ as a commuting set of derivations, i.e., the derivations $\delta_i: R\rightarrow R$ are defined for all $i$ and $\delta_i\delta_j=\delta_j\delta_i$ for all $i,j$.

Let $\Theta$ be the free commutative monoid on the set $\Delta = \{\delta_1,\ldots, \delta_m\}$  of
derivation operators. The elements
\bes
\theta = \delta_1^{i_1}\ldots \delta_m^{i_m}
\ees
of the monoid $\Theta$ are called {\em derivative}
operators. The {\em order} of $\theta$ is defined as $|\theta|=i_1+\ldots+i_m$.

Let $R$ be a differential ring. Denote by $R^e$ the free left $R$-module with a basis $\Theta$. Every element $u\in R^e$ can be uniquely written in the form
\bes
u=\sum_{\theta\in \Theta} r_{\theta}\theta
\ees
with a finite number of nonzero $r_{\theta}\in R$. It is well known \cite{KLMP} that there is the unique ring structure  on $R^e$ defined by the relations
\bes
\delta_i r=r \delta_i+ \delta_i(r)
\ees
for all $i$ and $r\in R$. Note that the ring $R^e$ is generated by $R$ and $\Delta$.
 Every left module over $R^e$ is called a {\em differential module} over $R$ \cite{KLMP}. For this reason we call $R^e$ the {\em universal enveloping} ring of $R$.

 Obviously, $R$ is a left $R^e$-module and every $I\subseteq R$ is a differential ideal of $R$ if and only if $I$ is an $R^e$-submodule of $R$.

Let $x^\Theta=\{x^\theta | \theta\in \Theta\}$ be a set of symbols enumerated by the elements of $\Theta$. Consider the polynomial algebra $R[x^\Theta]$ over $R$ generated by the set of (polynomially) independent variables $x^\Theta$. It is easy to check that the derivations $\delta_i$ can be uniquely extended to a derivation of $R[x^\Theta]$ by $\delta_i(x^\theta)=x^{\delta_i\theta}$. Denote this differential ring
by $R\{x\}$; it is called the {\em ring of differential polynomials} in $x$ over $R$.

By adding additional variables, we can obtain the differential ring $R\{x_1,x_2,\ldots,x_n\}$ of differential polynomials in $x_1,x_2,\ldots,x_n$ over $R$. Let $M$ be the free commutative monoid generated by all elements $x_i^{\theta}$, where $1\leq i\leq n$ and $\theta\in \Theta$. The elements of $M$ are called {\em differential monomials} of $R\{x_1,x_2,\ldots,x_n\}$. Every element
$a\in R\{x_1,x_2,\ldots,x_n\}$ can be uniquely written in the form
\bes
a=\sum_{m\in M} r_m m
\ees
with a finite number of nonzero $r_m\in R$.

An analog of the Hilbert Basis Theorem for differential algebras is known as the Ritt-Raudenbush Basis Theorem \cite{Kaplansky}. The condition of maximality for differential ideals does not hold for finitely generated differential rings. The Ritt-Raudenbush Basis Theorem claims that if $R$ is a radically Notherian ring containing the field of rational numbers $\mathbb{Q}$ then $R\{x\}$ is also radically Noetherian, i.e., every radically closed differential ideal is finitely generated.

Examples of equationally Noetherian differential rings are:

(i) Every differential field of characteristic zero is equationally Noetherian. This is a direct corollary of the Ritt-Raudenbush Basis Theorem and is one of the basics of differential algebraic geometry \cite[Chapter IV]{Kolchin}.

(ii) Every differential domain of characteristic zero is equationally Noetherian by Proposition \ref{p1}(1) since it is embeddable into a differential field.

In the remaining part of the paper we always assume that $R$ is an arbitrary associative and commutative ring with identity containing the field of rational numbers $\mathbb{Q}=\mathbb{Q}1\subseteq R$. Consider $R$ as a $\Delta$-ring with  $\delta_i(R)=0$ for all $i$. Then  $R\{x_1,x_2,\ldots,x_n\}$ becomes a $R$-algebra. Denote by $C_m=R[x_1,\ldots,x_m]$ the polynomial algebra over $R$ in the variables $x_1,\ldots,x_m$ considered as a differential algebra with derivations $\delta_i\equiv \frac{\partial}{\partial x_i}$ for all $1\leq i\leq m$.

\begin{pr}\label{p2}
Let $C=R\{y_1,y_2,...,y_n,\ldots\}$  be the differential polynomial algebra over $R$ in the variables $y_1,y_2,...,y_n,\ldots$. Let
$f$ be a nonzero element of $C$. Then there exists a homomorphism 
$\varphi:C\longrightarrow C_m$ of differential $R$-algebras such that $\varphi(f)\neq 0$.
\end{pr} 
\Proof Suppose that there exists a nonzero $f\in C$ such that $\varphi(f)= 0$ for any homomorphism $\varphi : C\longrightarrow C_m$ of differential $R$-algebras. Then $f=0$ is a nontrivial identity of the differential algebra $C_m$. Since $\mathbb{Q}\subseteq R$, we can linearize the identities \cite{KBKA}. Linearizing $f$ we can get a nonzero polylinear element $g\in C$ such that $g=0$ is an identity of $C_m$. Suppose that $g$ contains the variables $y_1,\ldots,y_n$ and is polylinear in these variables.

For any $\theta = \delta_1^{i_1}\ldots \delta_m^{i_m}\in \Theta$ we set $\a(\theta)=(i_1,\ldots,i_m)\in\mathbb{Z}_+^m$ and  
\bes
X(\theta) =\frac{1}{i_1! \ldots i_m!}x_1^{i_1}\ldots x_m^{i_m}\in C_m.
\ees
 If $\theta,\theta_1\in \Theta$ then we define set $\theta\leq \theta_1$ if  $\a(\theta)\leq \a(\theta_1)$ with respect to the lexicographical order in $\mathbb{Z}_+^m$. 
Notice that $X(\theta)^{\theta}=1$ and $X(\theta)^{\theta_1}=0$ for all $\theta,\theta_1\in \Theta$ if $\theta<\theta_1$. 

For any polylinear differential monomial
\bes
u=y_1^{\theta_1}y_2^{\theta_2}...y_n^{\theta_n}
\ees
of degree $n$ set $\b(u)=(\a(\theta_1),\ldots,\a(\theta_n))\in \mathbb{Z}_+^{mn}$. If $u$ and $v$ are two polylinear monomials then set $u\leq v$ if $\b(u)\leq \b(v)$ with respect to the lexicographical order in $\mathbb{Z}_+^{mn}$. 

The polylinear element $g$ of degree $n$ can be uniquely written in the unique form
\bes
g=\a y_1^{\theta_1}y_2^{\theta_2}...y_n^{\theta_n} + \sum \alpha_i m_i, 
\ees
where $0\neq \a\in R$ and $m_i$ is a polylinear differential monomial such that $y_1^{\theta_1}\ldots y_n^{\theta_n}<m_i$ for all $i$. 

Notice that 
\bes
X(\theta_1)^{\theta_1}\ldots X(\theta_n)^{\theta_n})=1 \ \mathrm{and} \ m_i(X(\theta_1),\ldots,X(\theta_n))=0 
\ees
for all $i$. Consequently, 
\bes
g(X(\theta_1),\ldots,X(\theta_n))=\a\neq 0, 
\ees
that is $g=0$ is not an identity of $C_m$.  $\Box$

\begin{co}\label{c4} 
The variety of all differential algebras ($\Delta$-algebras)  over $R$ is generated by $C_m=R[x_1,\ldots,x_m]$.
\end{co}
\Proof By Proposition \ref{p2}, the algebra $C_m$ does not satisfy any nontrivial differential identity. $\Box$

\begin{theor}\label{t2} Let $A$ be a differential $R$-algebra without nilpotent elements satisfying the ascending chain condition for radically closed differential ideals. Then $A$ is equationally Noetherian. 
\end{theor}
\Proof Notice that $A\{y_1,\ldots,y_n\}=A\sqcup R\{y_1,\ldots,y_n\}$ is the coproduct of $A$ and $R\{y_1,\ldots,y_n\}$ in the variety of differential $R$-algebras. Let $S$ be any subset of $A\{y_1,\ldots,y_n\}$ and 
\bes
Z(S)=\{(a_1,\ldots,a_n)\in A^n | s(a_1,\ldots,a_n)= 0 \ \mathrm{for \ all} \ s\in S\}. 
\ees

Let $I(S)$ be the set of all elements $f$ of $A\{y_1,\ldots,y_n\}$ such that $f(a_1,\ldots,a_n)= 0$ for all $(a_1,\ldots,a_n)\in Z(S)$. 
Obviously, $S\subseteq I(S)$ and $I(S)$ is an ideal of $A\{y_1,\ldots,y_n\}$. Note that 
\bes
f(a_1,\ldots,a_n)^{\delta_i}=f^{\delta_i}(a_1,\ldots,a_n)
\ees
for all $f\in A\{y_1,\ldots,y_n\}$ and for all $a_1,\ldots,a_n\in A$. Consequently, $I(S)$ is a differential ideal of $A\{y_1,\ldots,y_n\}$.  

If $f^k\in I(S)$ then $f^k(a_1,\ldots,a_n)=0$ for all $(a_1,\ldots,a_n)\in Z(S)$. The algebra $R$ does not contain nilpotent elements. Therefore  $f(a_1,\ldots,a_n)=0$ and $f\in I(S)$. Thus $I(S)$ is radically closed differential ideal of $A\{y_1,\ldots,y_n\}$. 

By the Ritt-Raudenbush basis Theorem \cite{Kaplansky, Kolchin, Ritt}, $A\{y_1,\ldots,y_n\}$ satisfies the ascending chain condition for radically closed differential ideals. Consequently, $I(S)$ coinsides with the radically closed differential ideal generated by a finite number of elements $f_1,\ldots,f_k\in A\{y_1,\ldots,y_n\}$. Therefore  $Z(S)=Z(f_1,\ldots,f_k)$. $\Box$

\begin{co}\label{c5} If $R$ is an assiciative and commutative radically Noetherian domain then the differential algebra  $C_m=R[x_1,\ldots,x_m]$ is equationally Noetherian. 
\end{co}
\Proof Notice that the differential structure on $R$ is trivial, i.e., $R$ is a radically Noetherian differential algebra. The algebra $C_m$ is radically Noetherian differential $R$-algebra by Ritt-Raudenbuch Basis Theorem \cite{Kaplansky, Kolchin, Ritt}. $\Box$

\begin{co}\label{c6} If $R$ is an associative and commutative radically Noetherian domain then every differential polynomial $R$-algebra is equationally Noetherian.  
\end{co}
\Proof Every free algebra of the variety of algebras $\mathrm{Var}(C_m)$ is equationally Noetherian by Proposition \ref{p1}. $\Box$

\section{The left-symmetric Witt algebra $L_m$ and the Witt algebra $W_m$}

\hspace*{\parindent}

{\bf Construction of left-symmetric and Lie algebras}. Consider differential algebras with the basic set of derivations $\Delta = \{\delta_1,\ldots, \delta_m\}$. Let $R$ be an arbitrary associative and commutative ring with identity and $A$ be an arbitrary differential algebra over $R$. Consider the space $B$ of elements of the form 
\bes
f_1\delta_1+f_2\delta_2+...+f_n\delta_m, \ \ f_i\in A, 1\leq i\leq m, 
\ees
of the universal enveloping algebra $C^e$. The product $\circ$ on $B$ defined by 
\bes
f\delta_i\circ g\delta_j=f\delta_i(g)\delta_j, \ \ \ f,g\in C,  1\leq i,j\leq m, 
\ees
turns $B$ into a left-symmetric algebra, that is $\langle B, \circ \rangle$ satisfies the identity 
 \bes
(x\circ y)\circ z-x\circ (y\circ z) = (y\circ x)\circ z -y\circ (x\circ z).
\ees
This can be easily checked as in \cite{KU16,16TG}. We denote the left-symmetric algebra $\langle B, \circ \rangle$ by $L(A)$ and call it the {\em left-symmetric Witt algebra of the differential algebra $A$}. 

The construction of $L(A)$ is a generalization  
 of the well known  Gelfand's construction  for Novikov algebras \cite{GD79}, and if $m=1$ then $L(A)$ is a Novikov algebra. 

Recall that the variety of left-symmetric algebras is Lie admissible, i.e., the commutator algebra of any left-symmetric algebra is a Lie algebra. Therefore  $B$ is a Lie algebra with respect to the bracket 
\bes
[f\delta_i, g\delta_j] = f\delta_i\circ g\delta_j -g\delta_j \circ f\delta_i. 
\ees
Denote the Lie algebra $\langle B, [\cdot,\cdot] \rangle$ by $W(A)$ and call it the {\em Witt algebra of the differential algebra $A$}. If $m=1$ then the bracket $[\cdot,\cdot]$ is called {\em the Wronskian} \cite{Poin18}. 

{\bf The left-symmetric Witt algebra $L_m$ and the Witt algebra $W_m$}. Let $C_m=R[x_1, x_2,\ldots, x_m]$ be the differential algebra with derivations $\delta_i= \partial_i$ for all $1\leq i\leq m$ on the polynomial algebra  $R[x_1, x_2,\ldots, x_m]$ . Then the left-symmetric algebra $L_m=L(C_m)$ is called the {\em left-symmetric Witt algebra of index $m$} and the Lie algebra $W_m=W(C_m)$, is the well known {\em Witt algebra of index $m$}  \cite{KU16,16TG}.  

{\bf Free algebras of the varieties $\mathrm{Var}(L_m)$ and $\mathrm{Var}(W_m)$}. 
Let $C=R\{y_{ij}| 1\leq i\leq n, 1\leq j\leq m\}$ be the differential polynomial algebra in the variables $y_{ij}$ where 
$1\leq i\leq n$ and $1\leq j\leq m$. 
We set 
\bes
Y_i=y_{i1}\delta_1+\ldots+y_{im}\delta_m, \ \ \ 1\leq i\leq n.
\ees
 Denote by ${\bf L}_n$ the subalgebra of the left-symmetric algebra $L(C)$ generated by $Y_1,\ldots,Y_n$, and by ${\bf W}_n$ the subalgebra of the Lie algebra $W(C)$  generated by $Y_1,\ldots,Y_n$.

Below we show that ${\bf L}_n$ is the free algebra of the variety of algebras $\mathrm{Var}(L_m)$ freely generated by $Y_1,\ldots,Y_n$ and ${\bf W}_n$  is the free algebra of the variety of algebras $\mathrm{Var}(W_m)$ freely generated by $Y_1,\ldots,Y_n$. 

{\bf "Coproducts".} Let $\widetilde{C}=C_m\{y_{ij}| 1\leq i\leq n, 1\leq j\leq m\}$ be the differential polynomial algebra over the differential algebra $C_m$ in the variables $y_{ij}$, where $1\leq i\leq n$ and $1\leq j\leq m$. Obviously, 
\bes
\widetilde{C}=C_m\otimes_R C
\ees
is the tensor product of the differential algebras $C_m$ and $C=R\{y_{ij}| 1\leq i\leq n, 1\leq j\leq m\}$. Moreover, these algebras can be identified with their images in $\widetilde{C}$ since they are free $R$-modules. Then the left-symmetric algebra $L(\widetilde{C})$ contains $L_m$ and ${\bf L}_n$, and the Lie algebra $W(\widetilde{C})$ contains $W_m$ and ${\bf W}_n$. 

 Denote by $\widetilde{L}_n$ the subalgebra of the left-symmetric algebra $L(\widetilde{C})$ generated by $L_m$ and  ${\bf L}_n$ and by $\widetilde{W}_n$ the subalgebra of the Lie algebra $W(\widetilde{C})$ generated by $W_m$ and ${\bf W}_n$.

\begin{pr}\label{p3} 
$(a)$ For any $F_1,\ldots,F_n\in L(\widetilde{C})$ there exists a unique homomorphism $\phi : \widetilde{L}_n\to L(\widetilde{C})$ of the left-symmetric algebras such that $\phi_{|L_m}=\mathrm{id}$ and $\phi(Y_i)=F_i$ for all $1\leq i\leq n$. 

$(b)$ For any $F_1,\ldots,F_n\in W(\widetilde{C})$ there exists a unique homomorphism $\phi : \widetilde{W}_n\to W(\widetilde{C})$ of the Lie algebras such that $\phi_{|W_m}=\mathrm{id}$ and $\phi(Y_i)=F_i$ for all $1\leq i\leq n$. 
\end{pr}
\Proof Consider  
\bes
F_i=f_{i1}\delta_1+\ldots+f_{im}\delta_m, \ \  f_{ij}\in \widetilde{C}, \ 1\leq i\leq n, \ 1\leq j\leq m. 
\ees
Obviously there exists a unique endomorphism  $\psi : \widetilde{C}\longrightarrow   \widetilde{C}$ of the differential $R$-algebra $\widetilde{C}$ such that $\psi_{|C_m}=\mathrm{id}$ and $\psi(y_{ij})=f_{ij}$ for all $i,j$. This endomorphism can be uniquely extended to an endomorphism of the enveloping algebra $\widetilde{C}^e$. We denote this extension of $\psi$ by the same symbol, i.e., $\psi : \widetilde{C}^e\to \widetilde{C}^e$ is an endomorphism of the associative algebra $\widetilde{C}^e$ such that $\psi(c)=c$ for all $c\in C_m$, $\psi(y_{ij})=f_{ij}$, and $\psi(\delta_j)=\delta_j$ for all $i,j$. It is clear that $\psi(L(\widetilde{C}))\subseteq L(\widetilde{C})$ and 
\bes
\psi_{|\widetilde{L}_n} : \widetilde{L}_n \to L(\widetilde{C})
\ees
is a homomorphism of the left-symmetric algebras. We set 
\bes
\phi=\psi_{|\widetilde{L}_n} : \widetilde{L}_n \to L(\widetilde{C}). 
\ees
 Then $\phi(Y_i)=F_i$ for all $i$ and $\phi_{|L_m}=\mathrm{id}$. Such a homomorphism is unique since $\widetilde{L}_n$ is generated by $L_m$ and by $Y_1,\ldots,Y_n$. 

Repeating the same arguments we can prove the statement $(b)$ of the theorem. $\Box$

\begin{theor}\label{t3} 
$(a)$ The left-symmetric algebra ${\bf L}_n$ is the free algebra of the variety of algebras $\mathrm{Var}(L_m)$ freely generated by $Y_1,\ldots,Y_n$.

$(b)$ The Lie algebra ${\bf W}_n$ is the free algebra of the variety of algebras $\mathrm{Var}(W_m)$ freely generated by $Y_1,\ldots,Y_n$.
\end{theor}
\Proof 
Let $M$ be the free algebra of the variety of algebras $\mathrm{Var}(L_m)$ freely generated by $z_1,\ldots,z_n$. 
It is well known that $M$ is isomorphic to a subalgebra of a direct power of $L_m$ \cite{Cohn81,Malcev73}. 

Let $I$ be an arbitrary set of indexes,  $L_m^I$ be the direct $I$-power of $L_m$, and assume that $M$ is isomorphic to a subalgebra of $L_m^I$. 
Let $F_1,\ldots,F_n$ be any set of elements of $L_m^I$. Suppose that $F_i=(f_{ij})_{j\in I}$. By Proposition \ref{p3} for each $j\in I$ there exists a homomorphism $\psi_j : \widetilde{L}_n\to L_m$ such that $\phi_j(Y_i)=f_{ij}$ for all $i$. For  $\phi_j={\phi_j}_{|{{\bf L}_n}}$ we have that 
\bes
\phi : {\bf L}_n\to L_m^I, 
\ees
defined by $\phi(a)=(\phi_j(a))_{j\in I}$ for all $a\in {\bf L}_n$, is a homomorphism of the left symmetric algebras such that $\phi(Y_i)=F_i$ for all $i$.

 It follows that there exists a homomorphism  
\bes
\tau : {\bf L}_n \to M, 
\ees
such that  $\tau(Y_i)=z_i$ for all $1\leq i\leq n$.

If $\tau$ is not an isomorphism then there exists a nonzero element $F=F(Y_1,\ldots,Y_n)$ such that $\phi(F)=0$ for all homomorphisms $\phi : {\bf L}_n\to L_m$. Consider  
\bes
F=f_1\delta_1+\ldots+f_m\delta_m, f_i\in C. 
\ees
Suppose that $f_i\neq 0$. By Proposition \ref{p2}, there exists a homomorphism $\psi : C\to C_m$ such that $\psi(f_i)\neq 0$. 
As above, we can extend this homomorphism to $\psi : C^e\to C_m^e$. If $\phi=\psi_{|L} : L\to L_m$, then   
\bes
\phi(F)=\psi(f_1)\delta_1+\ldots+\psi(f_i)\delta_i+\ldots+\psi(f_m)\delta_m\neq 0
\ees
since $\psi(f_i)\neq 0$. This contradiction proves that the kernel of $\tau$ is $0$, i.e., $\tau$ is an isomorphism. 

This proves the statement $(a)$ of the theorem. 
The statement $(b)$ can be proved similarly. $\Box$

\begin{theor}\label{t4} Let $R$ be an arbitrary radically Noetherian domain containing the field of rational numbers $\mathbb{Q}=\mathbb{Q}1\subseteq R$. Then 

$(a)$ the left-symmetric Witt algebra $L_m$ over $R$ is equationally Noetherian;   

$(b)$ the Witt algebra $W_m$ over $R$ is equationally Noetherian.  
\end{theor}
\Proof By Proposition \ref{p3}, for any natural $n\geq 0$ and for any system of elements $F_1,\ldots,F_n\in L(\widetilde{C})$  there exists a homomorphism $\phi : {\bf L}_n\to L(\widetilde{C})$ of the left-symmetric algebras such that  $\phi(Y_i)=F_i$ for all $1\leq i\leq n$. Therefore  $L(\widetilde{C})$ belongs to the variety of algebras $\mathrm{Var}(L_m)$.  

Let $L_m\sqcup {\bf L}_n$ be the coproduct of $L_m$ and ${\bf L}_n$ in the variety of algebras $\mathrm{Var}(L_m)$. 
Since  $\widetilde{L}_n$ is the subalgebra of $L(\widetilde{C})$ generated by $L_m$ and ${\bf L}_n$, there exists a unique homomorphism  
\bes
\phi : L_m\sqcup {\bf L}_n\to \widetilde{L}_n
\ees
such that $\phi_{|L_m}=\mathrm{id}$ and $\phi_{|{\bf L}_n}=\mathrm{id}$.

Let $S$ be any subset of $L_m\sqcup {\bf L}_n$ and 
\bes
Z(S)=\{(a_1,\ldots,a_n)\in L_m^n | s(a_1,\ldots,a_n)= 0 \ \mathrm{for \ all} \ s\in S\}. 
\ees
Recall that  $s(a_1,\ldots,a_n)$ is the image of $s$ under the homomorphism $L_m\sqcup {\bf L}_n\to L_m$ that is identical on 
$L_m$ and sends each $Y_i$ to $a_i$ for all $i$. By Proposition \ref{p3}, $s(a_1,\ldots,a_n)=0$ if and only if $\phi(s)(a_1,\ldots,a_n)=0$. Therefore 
\bes
Z(S)=\{(a_1,\ldots,a_n)\in L_m^n | \phi(s)(a_1,\ldots,a_n)= 0 \ \mathrm{for \ all} \ s\in S\}. 
\ees
Suppose that 
\bes
\phi(s)=s_1\delta_1+\ldots+s_m\delta_m, \ s_1,\ldots,s_m\in \widetilde{C},   
\ees
and 
\bes
a_i=a_{i1}\delta_1+\ldots+a_{im}\delta_m, \ a_{i1},\ldots,a_{im}\in C_m, \ 1\leq i\leq n. 
\ees
Then 
\bes
\phi(s)(a_1,\ldots,a_n)=\sum_{k=1}^m s_k(a_{ij})\delta_k. 
\ees

Let $I(S)$ be the set of all elements $f$ of $\widetilde{C}$ such that $f(a_1,\ldots,a_n)= 0$ for all $(a_1,\ldots,a_n)\in Z(S)$. 
Obviously, $s_i\in I(S)$ for all $s\in S$ and $1\leq i\leq m$. Moreover, $I(S)$ is an ideal of $\widetilde{C}$. We have 
\bes
f(a_{ij})^{\delta_k}=f^{\delta_k}(a_{ij})
\ees
for all $f\in \widetilde{C}$ and for all $a_{ij}\in C_m$. It follows that $I(S)$ is a differential ideal of $\widetilde{C}$.  

If $f^k\in I(S)$ then $f^k(a_{ij})=0$ for all $(a_1,\ldots,a_n)\in Z(S)$. Since $R$ is a domain, the algebras $C_m=R[x_1,\ldots,x_m]$ and $\widetilde{C}=C_m\{y_{ij}| 1\leq i\leq n, 1\leq j\leq m\}$ do not contain nilpotent elements. 
Therefore $f(a_{ij})=0$ and $f\in I(S)$. Thus $I(S)$ is radically closed differential ideal of $\widetilde{C}$. 

By the Ritt-Raudenbush Basis Theorem \cite{Kaplansky, Kolchin, Ritt}, $C_m$ and $\widetilde{C}$ satisfy the ascending chain condition for radically closed differential ideals. Therefore there exists a finite subset $T$ of $S$ such that the radically closed differential ideal of   $\widetilde{C}$ generated by all $t_i, 1\leq i\leq m, t\in T$ coinsides with $I(S)$. It  implies that 
 $Z(S)=Z(T)$. 

This proves the statement $(a)$ of the theorem. 
The statement $(b)$ can be proved similarly. $\Box$

\begin{co}\label{c7} Let $R$ be an arbitrary radically Noetherian domain containing the field of rational numbers $\mathbb{Q}=\mathbb{Q}1\subseteq R$. Then 

$(a)$ every free algebra over $R$  of the variety of $\mathrm{Var}(L_m)$  is equationally Noetherian;  

$(b)$ every free algebra over $R$  of the variety of $\mathrm{Var}(W_m)$  is equationally Noetherian.  

\end{co}
\Proof This is a direct corollary of Proposition \ref{p1} and Theorem \ref{t4}. $\Box$

\section{Symplectic Poisson algebra $P_m$}

\hspace*{\parindent}

{\bf Construction of Poisson algebras}. Consider differential algebras with the basic set of derivations $\Delta = \{\delta_1,\ldots, \delta_{2m}\}$ containing $2m$ elements.  Let $R$ be an arbitrary associative and commutative ring with identity and let $A$ be an arbitrary differential algebra over $R$. 
For any $f,g\in A$ we set 
\bes
\{f,g\}=\sum_{i=1}^m(\delta_{2i-1}(f)\delta_{2i}(g)-\delta_{2i}(f)\delta_{2i-1}(g)).
\ees

It is easy to show that the $R$-module $A$ with respect to the product $\cdot$ and the bracket $\{\cdot,\cdot\}$ is a Poisson algebra \cite{KMLS}, in other words,

$(a)$ $\langle A, \cdot\rangle$ is an associative and commutative algebra,

$(b)$ $\langle A, \{\cdot,\cdot\}\rangle$ is a Lie algebra, and 

$(c)$ $\langle A, \cdot,  \{\cdot,\cdot\}\rangle$ satisfies the (Leibniz) identity: 
\bes
\{xy,z\}=x\{y,z\}+\{x,y\}y. 
\ees

Denote the Poisson algebra $\langle A, \cdot, \{\cdot,\cdot\}\rangle$ by $P(A)$.  

{\bf The symplectic Poisson algebra $P_m$}. Let $C_{2m}=R[x_1,\ldots,x_{2m}]$ be the differential algebra with derivations 
 $\delta_i = \partial_i$ for all $1\leq i\leq 2m$ on the polynomial algebra $R[x_1,\ldots,x_{2m}]$. We denote  $P_m=P(C_{2m})$. The Poisson algebra $P_m$ is called the {\em symplectic Poisson algebra of index $m$ over $R$}.  If 
 $f,g\in P_m$ then
\bes
\{f,g\}=\sum_{i=1}^m( \frac{\partial f}{\partial x_{2i-1}}\frac{\partial g}{\partial x_{2i}}-\frac{\partial f}{\partial x_{2i}}\frac{\partial g}{\partial x_{2i-1}}).
\ees

{\bf Free algebras of the variety $\mathrm{Var}(P_m)$}. 
Let $C=R\{y_1,\ldots,y_n\}$ be the differential polynomial algebra over $R$ in the variables $y_1,\ldots,y_n$.  Denote by ${\bf P}_n$ the subalgebra of the Poisson algebra $P(C)$ generated by $y_1,\ldots,y_n$.  

Below we show that ${\bf P}_n$ is the free algebra of the variety of algebras $\mathrm{Var}(P_m)$ freely generated by $y_1,\ldots,y_n$.

{\bf "Coproducts".} Let $\widetilde{C}=C_{2m}\{y_1,\ldots,y_n\}$ be the differential polynomial algebra over the differential algebra $C_{2m}$ in the variables $y_1,\ldots,y_n$. It is clear that  
\bes
\widetilde{C}=C_{2m}\otimes_R C
\ees
is the tensor product of the differential algebras $C_{2m}$ and $C=R\{y_1,\ldots,y_n\}$. Moreover, these algebras can be identified with their images in $\widetilde{C}$ since they are free $R$-modules. Then the Poisson algebra $P(\widetilde{C})$ contains $P_m$ and ${\bf P}_n$.  

 Denote by $\widetilde{P}_n$ the subalgebra of the Poisson algebra $P(\widetilde{C})$ generated by $P_m$ and  ${\bf P}_n$.

\begin{pr}\label{p4} 
For any $F_1,\ldots,F_n\in P(\widetilde{C})$ there exists the unique homomorphism $\phi : \widetilde{P}_n\to P(\widetilde{C})$ of the Poisson algebras such that $\phi_{|P_m}=\mathrm{id}$ and $\phi(y_i)=F_i$ for all $1\leq i\leq n$. 
\end{pr}

The proof of this proposition repeats the proof of Proposition \ref{p3}. The proofs of the following two theorems are similar to the proofs of Theorem \ref{t3} and Theorem \ref{t4}, respectively.

\begin{theor}\label{t5} 
The Poisson algebra ${\bf P}_n$ is the free algebra of the variety of algebras $\mathrm{Var}(P_m)$ freely generated by $y_1,\ldots,y_n$.
\end{theor}

\begin{theor}\label{t6} Let $R$ be an arbitrary radically Noetherian domain containing the field of rational numbers $\mathbb{Q}=\mathbb{Q}1\subseteq R$. Then  the Poisson algebra $P_m$ over $R$ is equationally Noetherian. 
\end{theor}

\begin{co}\label{c8} Let $R$ be an arbitrary radically Noetherian domain containing the field of rational numbers $\mathbb{Q}=\mathbb{Q}1\subseteq R$. Then every free algebra over $R$  of the variety  $\mathrm{Var}(P_m)$  is equationally Noetherian.  
\end{co}

Using Proposition \ref{p3}, Theorem \ref{t3}, and general results on the structure of the variety of algebras generated by an algebra \cite{Cohn81,Malcev73}, it is possible to show that $\widetilde{L}_n$ is the coproduct of $L_m$ and ${\bf L}_n$ in the variety of algebras $\mathrm{Var}(L_m)$  and $\widetilde{W}_n$ is the coproduct of $W_m$ and ${\bf W}_n$ in the variety of algebras $\mathrm{Var}(W_m)$.  Similarly it can be proved that $\widetilde{P}_n$ is the coproduct of $P_m$ and ${\bf P}_n$ in the variety of algebras $\mathrm{Var}(P_m)$.

\bigskip

\begin{center}
	{\large Acknowledgments}
\end{center}

The first author is supported by Russian Science Foundation, grant \mbox{22-11-00052}.

The second author is supported by Nazarbayev University FDCRG (Grant No. 021220FD3851). 

The third author is supported by the Ministry of Education and Science of the Republic of Kazakhstan (project  AP14872073).

\end{document}